%% file: batson_counterexample.tex
\documentclass[10pt,reqno]{amsart}
\input{macros/preamble.tex}

\begin{document}
\thispagestyle{empty}
\title[A counterexample to Batson's conjecture]{A counterexample to Batson's conjecture.}
\author{Andrew Lobb} 
\address{Mathematical Sciences,
Durham University,
Durham,
UK.}
\email{andrew.lobb@durham.ac.uk}

\begin{abstract}
We show that the torus knot $T_{4,9}$ bounds a smooth M{\"o}bius band in the $4$-ball, giving a counterexample to Batson's non-orientable analogue of Milnor's conjecture on the smooth slice genera of torus knots.
\end{abstract}

\maketitle

Batson's conjecture says that the smooth non-orientable $4$-ball genus of a torus knot is realized by a simple construction.  This is analogous to Milnor's conjecture (verified by Kronheimer-Mrowka~\cite{km6}) that the smooth orientable $4$-ball genus of a torus knot is realized by the surface obtained from applying Seifert's algorithm to a standard diagram of the knot.

\section*{The conjecture.}
\label{sec:conjecture}

\input{sections/conjecture.tex}

\section*{A counterexample.}
\label{sec:countereg}
\input{sections/countereg.tex}

\bibliographystyle{amsplain}
\bibliography{references/works-cited.bib}
%
\end{document}

%% file: macros/preamble.tex
\usepackage{thmtools,enumerate,graphicx,xcolor}
\usepackage{etex,float}
\usepackage{ucs}
\usepackage{amsmath, amssymb, amscd, tabu, diagbox}
\usepackage{tikz}
\usetikzlibrary{arrows,chains,matrix,positioning,scopes}
\usepackage{mathtools}
\usepackage{sseq}
\usepackage{stmaryrd}
\usepackage[latin1]{inputenc}
\usepackage{amsmath}
\usepackage{amssymb}
\usepackage{epsfig}
\usepackage{psfrag}

\makeatletter
\def\@tocline#1#2#3#4#5#6#7{\relax
  \ifnum #1>\c@tocdepth 
  \else
    \par \addpenalty\@secpenalty\addvspace{#2}%
    \begingroup \hyphenpenalty\@M
    \@ifempty{#4}{%
      \@tempdima\csname r@tocindent\number#1\endcsname\relax
    }{%
      \@tempdima#4\relax
    }%
    \parindent\z@ \leftskip#3\relax \advance\leftskip\@tempdima\relax
    \rightskip\@pnumwidth plus4em \parfillskip-\@pnumwidth
    #5\leavevmode\hskip-\@tempdima
      \ifcase #1
       \or\or \hskip 2em \or \hskip 2em \else \hskip 3em \fi%
      #6\nobreak\relax
    \hfill\hbox to\@pnumwidth{\@tocpagenum{#7}}\par
    \nobreak
    \endgroup
  \fi}
\makeatother

\usepackage[all]{xy}

\usepackage{graphicx}
\usepackage{tabularx}
\usepackage{overpic}
\usepackage{tikz}
\usetikzlibrary{arrows}
\usepackage{rotating}
\usepackage{framed}
\usepackage{xcolor}
\usepackage{placeins}
\usepackage{booktabs}
\usepackage{pinlabel}
 \usepackage{picins}\usepackage{wrapfig}\usepackage{booktabs}\usepackage[font=footnotesize,labelfont=bf]{caption}
\usepackage{placeins}

\usepackage{color}
\usepackage{multirow}
\usepackage{overpic}
\usepackage{colortbl}
\usepackage[margin=1.4in]{geometry}

\makeatletter
\newcommand{\LeftEqNo}{\let\veqno\@@leqno}
\makeatother
\newlength{\lowerhalftmp}

\usepackage{microtype}

\theoremstyle{definition}

\renewcommand{\epsilon}{\varepsilon}

\DeclareMathAlphabet{\mathpzc}{OT1}{pzc}{m}{it}
\usepackage{mathrsfs}



%% file: sections/conjecture.tex
Let $T_{p,q} \subset S^3$ be the $(p,q)$ torus knot for $p>q\geq2$, and let $D_{p,q}$ be the usual $q$-stranded braid closure diagram of $T_{p,q}$.  Adding a blackboard-framed $1$-handle to two adjacent strands of $D_{p,q}$ results in a simpler torus knot, whose usual braid closure diagram we then consider.  Repeating this procedure eventually arrives at the unknot, which may be capped off in the~$4$-ball to give a surface $F_{p,q} \subset B^4$ with $\partial F_{p,q} = T_{p,q}$.  Batson conjectured \cite{batson} that $b_1(F_{p,q})$ is minimal among the first Betti numbers of non-orientable smooth surfaces in the $4$-ball with boundary~$T_{p,q}$.

We heard of this conjecture in a talk by Van Cott who, together with Jabuka, has verified it in many cases \cite{vcj}.

%% file: sections/countereg.tex
\begin{figure}[ht]
	\includegraphics[scale=0.035]{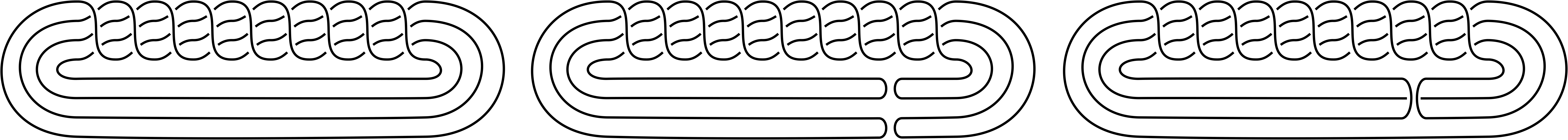}
	\caption{On the left is shown the torus knot $T_{4,9}$.  In the middle we have added two $1$-handles resulting in the unknot - this describes the surface $F_{4,9} \subset B^4$ which has $b_1(F_{4,9}) = 2$.  On the right we show how one may add a single $1$-handle to $T_{4,9}$ to result in the knot $6_1$, which is smoothly slice, thus giving a surface $\Sigma \subset B^4$ with $\partial \Sigma = T_{4,9}$ and $b_1(\Sigma) = 1$.}\label{fig:cobs}
\end{figure}